\documentclass{amsart}
\usepackage{amscd,amssymb,hyperref}
\usepackage[all]{xy}
\theoremstyle{plain}
\newtheorem{prop}[subsection]{Proposition}
\newtheorem{thm}[subsection]{Theorem}
\newtheorem{lem}[subsection]{Lemma}
\newtheorem{cor}[subsection]{Corollary}
\theoremstyle{remark}

\theoremstyle{definition}

\newcommand{\len}{len}
\numberwithin{equation}{section}

\title{On the degree $2$ map for a sphere }

\author[F.~R.~Cohen]{F.~R.~Cohen$^{*}$}
\address{Department of Mathematics, University of Rochester,
Rochester, NY 14627}
\email{\href{mailto:cohf@math.rochester.edu}{cohf@math.rochester.edu}}
\thanks{{$^{*}$}Partially supported by the National Science Foundation}

\author[I.~Johnson]{I.~Johnson}
\address{Department of Mathematics, Willamette University,
Willamette, OR 97301}
\email{\href{mailto:ijohnson@willamette.edu}{ijohnson@willamette.edu}}

\subjclass{Primary: , Secondary:}
\begin{document}
\begin{abstract}

The purpose of this article is to compare the two self-maps of
$\Omega^kS^{2n+1}$ given by $\Omega^k[2]$ the $k$-fold looping of a
degree $2$ map and $\Psi^k(2)$ the H-space squaring map. The main
results give that in case $2n+1 \neq 2^j-1$, these maps are
frequently not homotopic and also that their homotopy theoretic
fibres are not homotopy equivalent.

The methods are a computation of an unstable secondary operation
constructed by Brown and Peterson in the first case and the Nishida
relations in the second case.

One question left unanswered here is whether the maps
$\Omega^{2n+1}[2]$ and $\Psi^{2n+1}(2)$ are homotopic on the level
of $\Omega^{2n+1}_0S^{2n+1}$. A natural conjecture is that these two
maps are homotopic.
\end{abstract}

\maketitle

\section{Introduction and statement of results}

Consider the two natural self-maps of $\Omega^k S^{2n+1}$ given by
the $k$-fold looping of a degree $2$ map $$\Omega_{}^k[2]:\Omega^k
S^{2n+1} \to\
  \Omega^{k} S^{2n+1}$$ and
if $k \geq 1$, the $H$-space squaring map $$\Psi_{}^k(2):
  \Omega^k S^{2n+1} \to \Omega^{k} S^{2n+1}.$$ Furthermore,
  let $$\Omega^kS^{2n+1}\{[2]\}$$ denote the homotopy
theoretic fibre of $\Omega^k[2]$ and $$\Omega^kS^{2n+1}\{\Psi\}$$
denote the homotopy theoretic fibre of $\Psi^k(2)$. The purpose of
this article is to compare

\begin{itemize}
  \item the maps $\Omega_{}^k[2]$ and $\Psi_{}^k(2)$, as well as
  \item the spaces $\Omega^kS^{2n+1}\{\Psi\}$ and
  $\Omega^kS^{2n+1}\{[2]\}$.
\end{itemize}

The main point of this article is that the previous comparison is a
further illustration of the dichotomy between spheres of dimension
$2^j-1$ and spheres of other dimensions. Namely, Theorems
\ref{thm:theorem one}, and \ref{thm:theorem two} imply that if
$\Omega_{}^k[2]$, and $\Psi_{}^k(2)$ are homotopic, then the values
of $k$ are monotonically increasing with $n$ for certain choices of
$2n+1$ which are not equal to $2^j-1$. On the other hand, the strong
form of the Kervaire invariant conjecture implies that the maps
$\Omega_{}^k[2]$ and $\Psi_{}^k(2)$ are homotopic for $k=2$ in case
$2n+1 = 2^j-1$. Further discussion concerning this last point is
given in Proposition \ref{prop:prop divisibility of w} below in
which the dimensions of the spheres are $2n+1 = 2^j-1$.

\begin{thm}\label{thm:theorem one} Assume that the two self-maps
of $\Omega^k S^{2n+1}$ given by
$$\Omega_{}^k[2],\Psi_{}^k(2):\Omega^k S^{2n+1} \to \Omega^k
S^{2n+1}$$ are homotopic.
\begin{enumerate}
\item If $2n+1 = 2^{t+1}+ 2^t -1$ for $t \geq 1$, then $k  \geq 2^{t}+1$.
\item If $2n+1 = 2^{t+1}+1$ for $t \geq 1$, then $k \geq 2^{t} + 1$.
\end{enumerate}
\end{thm}

\begin{cor}\label{cor:corollary.of.theorem one}
There does not exist a finite integer $k$ such that for all $n>0$,
the two self-maps of $\Omega^k S^{2n+1}$ given by
$$\Omega_{}^k[2],\Psi_{}^k(2):\Omega^k S^{2n+1} \to \Omega^k
S^{2n+1}$$ are homotopic.
\end{cor}

Similar results are satisfied if the map $\Omega^k[2]$ is replaced
by $\Omega^k[-1]$ where $[-1]$ is a map of degree $-1$ and the map
$\Psi_{}^k(2)$ is replaced by $\Psi^k(-1)$, any choice of the loop
inverse for $\Omega^k S^{2n+1}$ \cite{1286}.

In addition, if the two self-maps of $\Omega^k S^{2n+1}$ given by
$\Omega_{}^k[2]$, and $\Psi_{}^k(2)$ are homotopic, then their
homotopy theoretic fibres are homotopy equivalent. However, the
converse does not appear to be evident. The next theorem, a direct
consequence of the Nishida relations, implies that these fibres are
frequently not homotopy equivalent in case $n$ is divisible by $8$.
Thus the maps $\Omega_{}^k[2]$ and $\Psi_{}^k(2)$ are not homotopic
in these cases. It seems likely that if the maps $\Omega_{}^k[2]$
and $\Psi_{}^k(2)$ are homotopic on $\Omega^k S^{2n+1}$ and $n$ is
restricted to values such that $2n+1 \neq 2^j-1$, then $\lim_{n \to
\infty } k = \infty$.

\begin{thm}\label{thm:theorem two}
If $n>1$ and $q\geq 1$, then the mod-$2$ cohomology of $\Omega^{2^n}
S^{q2^{n+2} + 1}\{\Psi\}$ and $\Omega^{2^n} S^{q2^{n+2} + 1}\{[2]\}$
are not isomorphic as modules over the mod-$2$ Steenrod algebra and
thus these spaces are not homotopy equivalent.
\end{thm}

\begin{cor}\label{cor:corollary.of.theorem two} If $n>1$ and
$q\geq 1$ then the maps
$$\Omega^{2^n}[2],\Psi^{2^n}(2):\Omega^{2^n} S^{q2^{n+2} + 1}
 \to \Omega^{2^n} S^{q2^{n+2} + 1}$$ are not homotopic.
\end{cor}

{\bf Remark:} Theorem \ref{thm:theorem one} sometimes gives stronger
information than Theorem \ref{thm:theorem two} concerning the
minimum values of $k$ such that the maps
$$\Omega^{k}[2],\Psi^{k}(2):\Omega^{k} S^{2^{n+2} + 1}
\to \Omega^{k} S^{2^{n+2} + 1}$$ are possibly homotopic. Both
results can be improved in special cases.

One example is given in Table $1$ below for the case of $\Omega^k
S^{17}$. The stated values of $k$ in Table $1$ are obtained as
consequences of the techniques used to prove Theorem
\ref{thm:theorem one} rather than the explicit statement of either
\ref{thm:theorem one}, or \ref{thm:theorem two}. In the case of
$\Omega^{k} S^{17}$, it follows that $k \geq 15$ by an application
of a secondary operation obtained from the relation
$$Sq^{18} = Sq^8[Sq^4(Sq^2Sq^4 + Sq^5Sq^1) + Sq^8Sq^2] + Sq^{16}Sq^2
+ Sq^{15}Sq^3+Sq^{14}Sq^4.$$ The requisite verifications are
sketched in section $2$ here after the proof of Theorem
\ref{thm:theorem one}.
\begin{table}[htdp] \caption{Results for $S^{2n+1}$,
$n=2, 4, 8, 5, 6$}
\begin{center}
\begin{tabular}{|c|}\hline \\
If $\Omega^k[2] \simeq \Psi^{k}(2) :  \Omega^kS^5 \rightarrow
\Omega^k S^5$, then $k \geq 3$. \\ \hline \\
If $\Omega^k[2] \simeq \Psi^{k}(2) : \Omega^k S^9 \rightarrow \Omega^k S^9$, then $k \geq 7$. \\ \hline \\
If $\Omega^k[2] \simeq \Psi^{k}(2) : \Omega^k S^{17} \rightarrow
\Omega^k S^{17}$, then $k \geq 15$. \\\hline \\
If $\Omega^k[2] \simeq  \Psi^{k}(2): \Omega^k S^{11} \rightarrow \Omega^k S^{11}$, then $k \geq 5$. \\ \hline \\
If $ \Omega^k[2] \simeq \Psi^{k}(2) : \Omega^k S^{13} \rightarrow \Omega^k S^{13}$, then $k \geq 7$. \\
\hline\end{tabular}\end{center}
\end{table}

The main results concerning these maps are closely tied to features
of the Whitehead square
$$[\iota_{n+t},\iota_{n+t}]: S^{2n+2t-1} \to\ S^{n+t}$$ denoted
$w_{n+t}$, as well as the classical James-Hopf invariant map
$$h_2: \Omega S^{n+t} \to \Omega S^{2n+ 2t -1}.$$ Most of the next result
is proven in \cite{1286} Proposition $11.3$ ( in which there is a
misprint where $\Omega^q(\phi)$ should be $\Omega^{q-1}(\phi)$ ). A
mildly different proof is included in section $4$ here for the
convenience of the reader.

\begin{prop}\label{prop:prop four}
If the two self-maps of $\Omega^k S^{2n+1}$ given by
$$\Omega_{}^k[2],\Psi_{}^k(2):\Omega^k S^{2n+1} \to
\Omega^k S^{2n+1}$$ are homotopic, then
\begin{enumerate}
\item $\Omega^k[\iota_{2n+1},\iota_{2n+1}]\circ
\Omega^{k-1}(h_2):\Omega^k S^{2n+1} \to\ \Omega^k S^{2n+1}$ is
null-homotopic and
    \item the composite denoted $\lambda(n,k)$
\[
\begin{CD}
\Sigma^{2n+1}(\mathbb R\mathbb P^{2n}/\mathbb R\mathbb P^{2n-k})=
\Sigma^{2n+1}(\mathbb R\mathbb P^{2n}_{2n-k+1})
@>{pinch}>> S^{4n+1}@>{w_{2n+1}}>> S^{2n+1}\\
\end{CD}
\] is null-homotopic.
\end{enumerate}

Furthermore, the difference $\Omega_{}^k[2]-\Psi_{}^k(2) = \Delta_k$
restricted to the second May-Milgram filtration of $\Omega^k
S^{2n+1}$ is null-homotopic if and only if the composite
$\lambda(n,k)$
\[
\begin{CD}
\Sigma^{2n+1}(\mathbb R\mathbb P^{2n}_{2n-k+1})
@>{pinch}>> S^{4n+1}@>{w_{2n+1}}>> S^{2n+1}\\
\end{CD}
\] is null-homotopic. Thus if $\lambda(n,k)$ is null-homotopic, there is a homotopy commutative diagram

\[
\begin{CD}
S^{4n+1}@>{}>>\Sigma^{2n+2}(\mathbb
R\mathbb P^{2n-1}_{2n-k+1})\\
 @V{w_{2n+1}}VV         @VV{\bar w_{2n+1}}V     \\
S^{2n+1} @>{1}>> S^{2n+1}.
\end{CD}
\] for some map $\bar w_{2n+1}$.

\end{prop}

It is reasonable to ask whether the two natural self-maps of
$\Omega^{2n+1}_0 S^{2n+1}$ given by
\begin{itemize}
  \item the $2n+1$-fold looping of a degree $2$ map, $\Omega^{2n+1}[2]$, and
  \item the $H$-space squaring map $\Psi^{2n+1}(2)$
\end{itemize} are homotopic, or whether the mod-$2$
cohomology algebras of $\Omega^{2n+1}_0S^{2n+1}\{\Psi\}$, and
$\Omega^{2n+1}_0S^{2n+1}\{[2]\}$ are isomorphic. The few known cases
occur for $2n +1 = 1,3,7,15,31,63$ as listed in \cite{1286,Selick}
with the case $2n+1 = 63$ based on the computations in
\cite{Kochman}. If these maps are homotopic, then the degree $2$ map
on $S^{2n+1}$ induces multiplication by $2$ on the level of homotopy
groups.

\begin{prop}\label{prop:prop divisibility of w}
Assume that $n \geq 0$.
\begin{enumerate}
\item The two self-maps of $\Omega S^{2n+1}$ given by $\Omega
[2]$, and $\Psi^1(2)$ are homotopic if and only if $w_{2n+1}= 0$.
Thus these two self-maps are homotopic if and only if $2n+1$ equals
$1$,$3$, or $7$.

\item The two self-maps of $\Omega^2 S^{2n+1}$ given by
$\Omega_{}^2[2]$, and $\Psi_{}^2(2)$ are homotopic if and only if
the Whitehead square $w_{2n+1} = [\iota_{2n+1}, \iota_{2n+1}]$ is
divisible by $2$. Thus these maps
\begin{enumerate}
    \item are homotopic for $n = 1,3,7,15,31,63$ and
    \item are not homotopic when $2k+1$ is not $2^j -1$ for some $j$.
\end{enumerate}
\end{enumerate}
\end{prop}

Further information concerning the divisibility of the Whitehead
square is listed next. In case $2n+1 = 2^j-1$, that the Whitehead
square is divisible by $2$ is known as the strong form of the
Kervaire invariant one conjecture and is known to admit a positive
solution in case $2n+1$ is $1$, $3$, $7$, $15$, $31$, or $63$
\cite{Kochman, Mahowald, Selick}. The Whitehead square is not
divisible by $2$ in case $2n+1 \neq 2^j-1$. Little is known about
the answer in general in case $2n+1 = 2^j-1 > 63$. A reformulation
of this topology question in terms of the Lie group $G_2$ and the
zero divisors in a classical construction of L.~E.~Dickson
concerning a (usually non-associative) multiplication on $\mathbb
R^{2^n}$ is given in \cite{Adem proceedings} with additional
information given in \cite{Moreno}.

In view of these examples, it appears that the cohomology  of
$\Omega^{2n+1}_0 S^{2n + 1}\{\Psi\}$ and $\Omega^{2n +1 }_0S^{2n +
1}\{[2]\}$ might be interesting as algebras over the Steenrod
algebra. The homology of $\Omega^{2n + 2 }_0S^{2n + 1}$ has been
worked out by T.~Hunter \cite{Hunter}.

The proof of Theorem \ref{thm:theorem one} depends on factorizations
of $Sq^{2n+2}$ for $2n+2 \neq 2^j$ which give lower bounds on $k$
via the method of evaluation of secondary operations of Brown, and
Peterson \cite{BP}. For example, $Sq^{10}$ can be factored in the
following two ways
$$Sq^{10} = Sq^4Sq^2Sq^4+Sq^8Sq^2+Sq^4Sq^5Sq^1,
\textrm{ \hspace{.2in} and \hspace{.2in} } Sq^{10} =
Sq^2Sq^8+Sq^9Sq^1.$$  An application of the Brown and Peterson
secondary operation with the first factorization gives that
$\Omega^6[2] \not\simeq \Psi^6(2): \Omega^6S^9 \rightarrow
\Omega^6S^9$, while an application of the operation with the second
factorization gives that $\Omega^2[2]\not\simeq \Psi^2(2): \Omega^2
S^9 \rightarrow \Omega^2 S^9$. Hence the first factorization gives a
stronger result.

The best lower bounds for $k$ using the methods above occur for the
smallest value of $k$ such that $Sq^{2n+2}$ is in the left ideal of
the Steenrod algebra given by $$L(k) = \mathcal{A} \{ Sq^1, Sq^2,
\dots, Sq^{2^k} \}.$$ In the example above, $Sq^{10}$ is in $L(2)$,
and in fact $k=2$ is the smallest such $k$. The smallest value of
$k$ such that $Sq^{2n+2}$ is in $L(k)$ is given in \cite{JM}, and is
described next.

\begin{thm}\label{thm:johnson.merzel}
If the two self-maps of $\Omega^k S^{2n+1}$ given by $\Omega^k[2]$,
and $\Psi^k(2)$ are homotopic, then $$k \geq F(2n+2),$$ where the
function $F$ is defined below.
\end{thm}

The following notation is used to define the function $F$. Given any
positive integer $n$, let $[n]$ denote the dyadic expansion of $n$
viewed as an ordered sequence of zeroes and ones with right
lexicographical ordering. That is if $$n = 2^{j_t}+2^{j_{t-1}} +
\cdots + 2^{j_1}+2^{j_0}$$ with $$j_t > j_{t-1} > \cdots > j_1
> j_0 \geq 0,$$ then the dyadic expansion of $n$ is ambiguously
denoted $$\alpha= (\alpha_q,\alpha_{q-1},\cdots,\alpha_1,\alpha_0)$$
for which
\[
\alpha_s=
\begin{cases}
1 & \text{if $s = j_m$ for $t \geq m \geq 0$, and}\\
0 & \text{if $s \neq j_m $ for $t \geq m \geq 0$.}
\end{cases}
\] Notice that $(\alpha_q,\alpha_{q-1},\cdots,\alpha_1,\alpha_0)$ and
$(0,\alpha_q,\alpha_{q-1},\cdots,\alpha_1,\alpha_0)$ are dyadic
expansions of the same integer, and are both equal to $[n]$.

Given a binary string $\alpha$,
\begin{enumerate}
\item let $|\alpha|$ denote the integer with binary expansion
$\alpha$,

\item let $\len(\alpha)$ denote the length of the binary string
$\alpha$, and \item let $z(\alpha)$ denote the number of
non-trailing zeroes in a string $\alpha$, thus if $$\alpha=
(\alpha_q,\alpha_{q-1},\cdots,\alpha_1,\alpha_0)$$ as above, then
$z(\alpha)= \len(\alpha)-t-j_0$. ( For example, $$|010010000|= 2^7 +
2^4 = |10010000|,$$  $$\len(010010000) = 9,$$ and $z(010010000)$ is
$3$.)
\end{enumerate}

Given binary strings $\alpha$ and $\beta$, let $\alpha \beta$ denote
their concatenation.  With this notation the function $F$ is defined
on a positive integer $n$ as follows. Write $[n]=\alpha \beta$ such
that $|\alpha| < z(\beta)$ and $\len(\beta)$ is minimal. Then
$$F(n)=n-2^{\len(\beta)-2}+1.$$

There are two main computations in this article. One is the
evaluation of an unstable secondary operation due to Brown and
Peterson \cite{BP} which gives a proof of Theorem \ref{thm:theorem
one}. The second is a computation with the Nishida relations
\cite{Nishida} which gives a proof of Theorem \ref{thm:theorem two}.

A table of contents of this paper is as follows:
\begin{enumerate}
\item [1:] Introduction and statement of results

\item [2:] Unstable secondary operations related to the Whitehead
product and the proof of Theorem \ref{thm:theorem one}

\item [3:]The Nishida relations and
the proof of Theorem \ref{thm:theorem two}

\item [4:]On the Proof of Proposition \ref{prop:prop four}

\end{enumerate}

The authors would like to thank J\'esus Gonzalez, Miguel
Xicot\'encatl as well as the other organizers for an interesting and
fruitful conference. The authors are grateful to the referee for his
careful reading of this article as well as for numerous excellent
suggestions. Finally, the authors would like to congratulate Sam
Gitler on his birthday and to thank him for his interest and
contagious joy in doing mathematics.

\section{Unstable secondary operations related to the Whitehead product, and the proof
of Theorem \ref{thm:theorem one}}

To prove Theorem \ref{thm:theorem one}, it suffices to check the
statement that if the two self-maps of $\Omega^k S^{2n+1}$ given by
$\Omega_{}^k[2],\Psi_{}^k(2):\Omega^k S^{2n+1} \to \Omega^k
S^{2n+1}$ are homotopic, then

\begin{itemize}
    \item for the case $2n+1 = 2^{t+1}+ 2^t -1$ with $t \geq 1$,
    it follows that $k \geq 2^{t}+1$ and
\item for the case $2n+1 = 2^{t+1}+1$ for $t \geq 1$, it follows that $k \geq 2^t + 1$.
\end{itemize} The steps in the strategy of the proof are as follows.
\begin{enumerate}
  \item Consider the difference of the two maps $\Omega_{}^k[2]$ and
$\Psi_{}^k(2)$ restricted to the second May-Milgram filtration of
$\Omega^k S^{2n+1}$.
  \item Observe that this difference factors through the suspension
  of a truncated projective space.
    \item Compute a non-trivial unstable secondary operation in the
    cohomology of the suspension of a truncated projective space
  (for the special cases listed directly above).
  \item Conclude that the difference is essential when restricted to
  the suspension of a truncated projective space as given in step $3$.
  \item Conclude by Proposition \ref{prop:prop four} that the two maps fail to be homotopic
  when restricted to the second May-Milgram filtration of
$\Omega^k S^{2n+1}$ in these special cases.
\end{enumerate} The proof of Theorem \ref{thm:theorem one} thus reduces to an evaluation
of certain unstable secondary cohomology operations constructed by
Brown, and Peterson \cite{BP}.

The Adem relations are $Sq^{i}Sq^{j} = \Sigma_{0 \leq s \leq [i/2]}
\left({}^{j-s-1}_{i-2s} \right) Sq^{i+j - s}Sq^{s}$ in case $i<2j$.
Thus if $2n+2$ is not a power of $2$, $Sq^{2n+2}$ appears in some
Adem relation which can be rewritten as
$$Sq^{2n+2} = \sum_{t_i \neq 2n+2} a_iSq^{t_i}$$ for which $a_i$ is
in the Steenrod algebra. The special cases given below suffice for
the purposes here with more complete answers given in \cite{JM}.

The secondary cohomology operations, devised by Brown, and Peterson
\cite{BP} to detect the Whitehead square on spheres not of dimension
$2^k-1$ are described next in order to set up the context for the
applications here. The results of \cite{BP} also give tertiary
operations which detect the Whitehead square for spheres of
dimension $2^k-1$ with $ k > 3$, but these operations are not used
here. Additional information concerning secondary operations is
listed in \cite{Harper}.

Consider the Eilenberg-Mac Lane space $K(\mathbb Z/2\mathbb Z,n)$
together with factorizations of $Sq^{2n+2} = \sum_{t_i \neq 2n+2}
a_iSq^{t_i}$ to obtain
\[
\begin{CD}
K(\mathbb Z/ 2 \mathbb Z, 2n+2)  @>{\Pi_{t_i \neq 2n+2
}Sq^{{t_i}}}>> \Pi_{t_i \neq 2n+2}K(\mathbb Z/2\mathbb Z, 2n+2+t_i)
\end{CD}
\] with homotopy theoretic fibre denoted ambiguously by $E_{2n+2}$
( depending on the choice of factorization of $Sq^{2n+2}$ ). Thus
there is a fibration $$\Pi_{t_i \neq 2n+2}K(\mathbb Z/2\mathbb
Z,2n+1+t_i) \to\ E_{2n+2} \to\ K(\mathbb Z/ 2 \mathbb Z, 2n+2) $$
for which $\iota_{2n+2}$ denotes the fundamental cycle for the base,
and $x_{2n+1+t_i}$ denotes the fundamental cycle for each
Eilenberg-Mac Lane space in the fibre. The transgression of the
cohomology class $\Sigma_i a_i x_{2n+1+t_i}$ is
$$\Sigma_i a_i Sq^{t_i}(\iota_{2n+2}) = Sq^{2n+2}(\iota_{2n+2})=
\iota_{2n+2}^2.$$

Next, consider the looped fibration above to obtain an analogous
fibration $$ \Pi_{t_i \neq 2n+2}K(\mathbb Z/2\mathbb Z,2n+t_i) \to\
\Omega E_{2n+2} \to\ K(\mathbb Z/ 2 \mathbb Z, 2n+1).$$ The
analogous cohomology class $$\Sigma_i a_i x_{2n+t_i}$$ obtained for
this last fibration is an infinite cycle in the Serre spectral
sequence for this fibration. Thus there is a choice of cohomology
class $\Phi(\iota_{2n+1})$ in the cohomology of $\Omega E_{2n+2}$
which restricts to $\Sigma_i a_i x_{2n+t_i}$ in the cohomology of
$\Pi_{t_i \neq 2n+2}K(\mathbb Z/2\mathbb Z,2n+t_i)$. Brown, and
Peterson \cite{BP} show that the reduced coproduct for
$\Phi(\iota_{2n+1})$ is non-trivial, and is thus given by
$$\iota_{2n+1} \otimes \iota_{2n+1}$$ by degree considerations.

The first non-vanishing homotopy group of $E_{2n+2}$ is given by
$\pi_{2n+2}(E_{2n+2}) = \mathbb Z/ 2 \mathbb Z$ with a choice of
representative for a generator denoted
$$\lambda: S^{2n+2} \to E_{2n+2}.$$ The loop map
$$f = \Omega(\lambda): \Omega S^{2n+2} \to\ \Omega E_{2n+2}$$
induces an isomorphism $f_*:H_{2n+1}(\Omega S^{2n+2}) \to
H_{2n+1}(\Omega E_{2n+2})$.

Let ${b_{2n+1}}_*$ denote the fundamental cycle in $H_{2n+1}(\Omega
S^{2n+2})$. Recall that $({b_{2n+1}}_*)^2$ is non-zero in the
Pontrjagin ring. Thus $f^*(\Phi(\iota_{2n+1}))$ evaluates
non-trivially on $({b_{2n+1}}_*)^2$ with the natural pairing
$<f^*(\Phi(\iota_{2n+1})), ({b_{2n+1}}_*)^2 > = 1$. Consequently,
$f$ induces a non-trivial map $$f_*:H_{4n+2}(\Omega S^{2n+2}) \to
H_{4n+2}(\Omega E_{2n+2}).$$

Furthermore, since the indeterminacy of the operation is given by
any choice of map to the fibre $$\Omega S^{2n+2} \to\ \Pi_{t_i \neq
2n+2}K(\mathbb Z/2\mathbb Z,2n+t_i),$$ the indeterminacy is always
trivial. Since the attaching map for the $4n+2$ cell in a minimal
cell decomposition of $\Omega S^{2n+2}$ is the Whitehead product
$w_{2n+1}$, the operation of \cite{BP} detects this element as long
as $2n+2 \neq 2^s$. These operations will be applied to the
following context.

Recall $F_s$ the $s$-th May-Milgram filtration \cite{May} of $\Omega
^{k}S^{2n+1}$  which was exploited earlier by Toda \cite{Toda} in
the special case of $s = 2$. The inclusion $F_{s-1}$ in $F_s$ is a
cofibration. The filtration quotient
$$F_2/F_1$$ is homotopy equivalent to $$D_2(\Omega ^{k}
S^{2n+1}) = \Sigma^{2n+1-k}(\mathbb R \mathbb P^{2n} / \mathbb R
\mathbb P^{2n-k}) = \Sigma^{2n+1-k}(\mathbb R \mathbb
P^{2n}_{2n+1-k}).$$

Next consider the cofibre sequence
\[
\begin{CD}
\Sigma^{2n+1-k}(\mathbb R\mathbb P^{2n-1}_{2n+1-k})@>{inclusion}>>
\Sigma^{2n+1-k}(\mathbb R\mathbb P^{2n}_{2n+1-k})@>{K}>> S^{4n+1-k}
\end{CD}
\] for which $$K:\Sigma^{2n+1-k}(\mathbb R\mathbb P^{2n}_{2n+1-k}) \to\
S^{4n+1-k}$$ denotes the natural collapse map with an induced
isomorphism (in mod $2$ homology) $$K_*:
H_{4n+1-k}(\Sigma^{2n+1-k}(\mathbb R\mathbb P^{2n}_{2n+1-k}))\to
H_{4n+1-k}(S^{4n+1-k}).$$ In addition, there is a ``boundary" map
obtained from the Barratt-Puppe sequence
$$\delta: S^{4n+1-k} \to\ \Sigma^{2n+2-k}(\mathbb R\mathbb P^{2n-1}_{2n+1-k}).$$

Consider the self-map of $\Omega ^{k}S^{2n+1}$ given by the
difference $$\Delta_{k} = \Omega^{k}[2] -\Psi^{k}(2)$$ restricted to
the second filtration $F_2 =F_2(\Omega ^{k} S^{2n+1})$. By
Proposition \ref{prop:prop four}, the difference $\Delta_{k} =
\Omega_{}^k[2]-\Psi_{}^k(2)$ restricted to the second May-Milgram
filtration of $\Omega^k S^{2n+1}$ is null-homotopic if and only if
the natural composite
\[
\begin{CD}
\Sigma^{2n+1}(\mathbb R\mathbb P^{2n}_{2n-k+1})
@>{K}>> S^{4n+1}@>{w_{2n+1}}>> S^{2n+1}\\
\end{CD}
\] denoted $\lambda(n,k)$ is null-homotopic in which case, there is a homotopy commutative diagram

\[
\begin{CD}
S^{4n+1}@>{}>>\Sigma^{2n+2}(\mathbb
R\mathbb P^{2n-1}_{2n - k +1})\\
 @VV{w_{2n+1}}V         @VV{\bar w_{2n+1}}V     \\
S^{2n+1} @>{1}>> S^{2n+1}.
\end{CD}
\] for some map $\bar w_{2n+1}$.

The next step in the proof of Theorem \ref{thm:theorem one} is an
examination of certain values of $n$ and $k$ such that the composite
$\lambda(n,k)$ is essential. By definition, there is a homotopy
commutative diagram
\[
\begin{CD}
\Sigma^{2n+1}(\mathbb R\mathbb P^{2n}_{2n-k+1})@>{\lambda(n,k)}>>
S^{2n+1}\\
 @V{K}VV         @VV{1}V     \\
S^{4n+1} @>{w_{2n+1}}>> S^{2n+1}
\end{CD}
\] together with an induced morphism of cofibre sequences

\[
\begin{CD}
\Sigma^{2n+1}(\mathbb R\mathbb P^{2n}_{2n-k+1})@>{\lambda(n,k)}>>
S^{2n+1} @>{}>>   X(n,k)    \\
 @V{K}VV         @VV{1}V     @VV{g}V  \\
S^{4n+1} @>{w_{2n+1}}>> S^{2n+1} @>{}>> J_2(S^{2n+1})
\end{CD}
\] where $J_2(S^{2n+1})$ denotes the second stage of the James
construction. Furthermore, any choice of map $$g:X(n,k) \to
J_2(S^{2n+1})$$ obtained from a morphism of cofibre sequences
induces an isomorphism $$g_*:H_{i}X(n,k) \to H_{i}J_2(S^{2n+1})$$
for $i = 2n+1, 4n+2$ by inspection of the long exact sequence in
homology obtained from a cofibre sequence. (Note: The choice of $g$
is not necessarily unique up to homotopy.)

Furthermore, there is a homotopy equivalence  after one suspension,
$$\Sigma( S^{2n+1} \vee \Sigma^{2n+2}(\mathbb R\mathbb
P^{2n}_{2n-k+1})) \to \Sigma X(n,k).$$ Thus the action of the
Steenrod algebra for the cohomology of $\Sigma X(n,k)$ is obtained
from the action on the cohomology of a truncated projective space.

Consider the composite
\[
\begin{CD}
X(n,k) @>{g}>> J_2(S^{2n+1}) @>{i}>>   \Omega S^{2n+2}\\
\end{CD}
\] denoted $$G: X(n,k) \to \Omega S^{2n+2}$$ for which $i$ is an equivalence through the $6n+2$ skeleton.
Recall the map $$f: \Omega S^{2n+2} \to\ \Omega E_{2n+2}$$ defined
earlier in this proof. Since
$$<f^*(\Phi(\iota_{2n+1})), ({b_{2n+1}}_*)^2> = 1,$$
it follows that $$<(G \circ f)^*(\Phi(\iota_{2n+1})), e_{4n+2}> =
1$$ where $e_{4n+2}$ is the unique non-trivial class in
$H_{4n+2}(X(n,k))$. Thus, if the indeterminacy of the choice of map
$G \circ f: X(n,k) \to \Omega E_{2n+2}$ is zero, then $\lambda(n,k)$
is essential as the resulting cohomology operation is itself
non-trivial. Vanishing of the indeterminacy will be checked next for
special cases.

Recall the Adem relations $Sq^{i}Sq^{j} = \Sigma_{0 \leq s \leq
[i/2]} \left({}^{j-s-1}_{i-2s} \right) Sq^{i+j - s}Sq^{s}$. There
are two such relations to consider in the proof of Theorem
\ref{thm:theorem one}. The first case is handled by the choice of
operation arising from the relation $$Sq^{2^t}Sq^{2^{t+1}} =
\Sigma_{0 \leq s \leq 2^{t-1}}\left({}^{2^{t+1}-s-1}_{2^t-2s}
\right) Sq^{2^t + 2^{t+1}- s}Sq^{s}$$ for $  t\geq 1$. It follows
that $$ Sq^{2^{t+1}+ 2^t } = Sq^{2^t}Sq^{2^{t+1}} + \Sigma_{1\leq s
\leq 2^{t-1}}\left({}^{2^{t+1}-s-1}_{2^t-2s} \right) Sq^{2^t +
2^{t+1}- s}Sq^{s}$$ with
\begin{enumerate}
\item $2n+2 = 2^{t+1}+ 2^t$ for $ t \geq 1$,
    \item $2^{t}-2s \geq 0$, and
    \item $2^t + 2^{t+1}- s \geq 2^t + 2^{t+1}- 2^{t-1} = 2^{t-1} + 2^{t+1}$.
\end{enumerate} To check that this operation has zero indeterminacy and thus
that the map $\lambda(n,k)$ is essential, it suffices to check that
the operation $a_s =Sq^{2^t + 2^{t+1}- s}$ vanishes on the
cohomology of $X(n,k)$ for $2n+2 = 2^{t+1}+ 2^t$, $ t \geq 1$ and
$2^{t}-2s \geq 0$ with $ s > 0$.

Since there is a homotopy equivalence $$\Sigma( S^{2n+1} \vee
\Sigma^{2n+2}(\mathbb R\mathbb P^{2n}_{2n-k+1})) \to \Sigma
X(n,k),$$ it suffices to check that the operation $a_s =Sq^{2^t +
2^{t+1}- s}$ vanishes on the cohomology of $\mathbb R\mathbb
P^{2n}_{2n-k+1}$. The assumption that $k-1 < 2^{t}$ gives that this
operation has zero indeterminacy by a check of degrees and thus the
map $\lambda(n,k)$ is essential.

The assumption that the two self-maps of $\Omega^k S^{2n+1}$ given
by $\Omega_{}^k[2]$, and $\Psi_{}^k(2)$ are homotopic implies that
$\Delta_{k}$, and hence $\lambda(n,k)$ is null in case $n = 2^t +
2^{t-1}-1$, with $k -1 <  2^{t} $, contradicting Proposition
\ref{prop:prop four}. Hence, $k-1 \geq 2^{t}$ and Theorem
\ref{thm:theorem one}, part $1$, follows.

The second case is handled by the choice of operation arising from
the relation   $$Sq^{2^t}Sq^{2+ 2^{t}} =\Sigma_{0 \leq s \leq
2^{t-1}} \left({}^{2^t+1-s}_{2^t-2s} \right) Sq^{2+ 2^{t+1} -
s}Sq^{s}$$ for $  t\geq 1$. It follows that $$Sq^{2+ 2^{t+1}}=
Sq^{2^t}Sq^{2+ 2^{t}} + \Sigma_{1 \leq s \leq 2^{t-1}}
\left({}^{2^t+1-s}_{2^t-2s} \right) Sq^{2+2^{t+1}- s}Sq^{s}$$ with

\begin{enumerate}
\item $2n + 2 = 2 + 2^{t+1}$ for $ t \geq 1$,
    \item $2^t \geq 2s$, and
    \item $2+2^{t+1}-s \geq 2+2^{t+1}- 2^{t-1} = 2 + 2^t +
    2^{t-1}.$
\end{enumerate} Thus if $k < 1 + 2^t $ with $n = 2^t$ and $2^{t-1} \geq s$ as
above, the operation $a_s = Sq^{2+2^{t+1}- s}$ vanishes on the
cohomology of $\mathbb R\mathbb P^{2n}_{2n-k+1}$ and thus $a_s$
vanishes on the cohomology of $X(n,k)$. Hence if $k < 1 + 2^t$, the
associated operation has zero indeterminacy in the cohomology of
$X(n,k)$ and so the map $\lambda(n,k)$ is essential. The rest of the
proof in this case is analogous to that for the first case and is
omitted. It follows that $k \geq 1 + 2^{t}$ and Theorem
\ref{thm:theorem one}, part $2$, follows.

The values of $k$ given in Table $1$ follow from an analogous
secondary operation obtained from an iteration of the Adem relations
as listed next rather than the explicit estimates in Theorem
\ref{thm:theorem one}. For example, the relation $$Sq^{18} =
Sq^8[Sq^4(Sq^2Sq^4 + Sq^5Sq^1) + Sq^8Sq^2] + Sq^{16}Sq^2 +
Sq^{15}Sq^3+Sq^{14}Sq^4$$ gives a stronger result than that stated
in Theorem \ref{thm:theorem one}. The indeterminacy of the
associated secondary operation is zero by an inspection of the
action of the Steenrod operations on the cohomology of $X(8,14)$ a
space which satisfies the property that $\Sigma X(8,14)$ is homotopy
equivalent to $\Sigma( S^{17} \vee \Sigma^{18}(\mathbb R\mathbb
P^{16}_{4}))$. The relations listed next are used to give the
results in Table $1$ by a direct check that the indeterminacy
vanishes in these cases. The details are omitted.
\begin{enumerate}
        \item $Sq^6 = Sq^2Sq^4 + Sq^5Sq^1$.
        \item $Sq^{10} = Sq^4Sq^6 + Sq^8Sq^2 = Sq^4(Sq^2Sq^4 + Sq^5Sq^1) + Sq^8Sq^2 $.
        \item $Sq^{18} = Sq^8Sq^{10} + Sq^{16}Sq^2 + Sq^{15}Sq^3+Sq^{14}Sq^4$ and thus
$Sq^{18} = Sq^8[Sq^4(Sq^2Sq^4 + Sq^5Sq^1) + Sq^8Sq^2] + Sq^{16}Sq^2
+ Sq^{15}Sq^3+Sq^{14}Sq^4$.
        \item $Sq^{12}=Sq^4Sq^8+ Sq^{11}Sq^1+ Sq^{10}Sq^2$.
        \item $Sq^{14} = Sq^6Sq^{8} + Sq^{13}Sq^1+Sq^{11}Sq^3$.
\end{enumerate}

\section{The Nishida relations and the proof of Theorem \ref{thm:theorem two}}

Information concerning the mod-$2$ homology of the spaces
$\Omega^kS^{n}\{[2]\}$ and $\Omega^kS^{n}\{\Psi\}$ is given below.
This information is used to show that the action of the Steenrod
operations on the mod-$2$ cohomology of the spaces in Theorem
\ref{thm:theorem two} differ, thus proving the Theorem. References
are \cite{CCPS, 1286}. In what follows below, assume that $1 < k <
n-3$ with homology always taken with coefficients in $\mathbb Z /2
\mathbb Z$.

The mod-$2$ homology $H_*\Omega^k S^n\{[2]\}$ is isomorphic to
$$H_*\Omega^k S^n \otimes H_*\Omega^{k+1}S^n$$ as a Hopf algebra
in case $1 < k  < n-3 $ with the natural map $$H_*\Omega^k
S^n\{[2]\}\to\ H_*\Omega^k S^n$$ induced by a $k$-fold loop map.
This map is an epimorphism of Hopf algebras over the mod-$2$
Steenrod algebra. Thus there is a unique class $x_{n-k}$ in
$H_{n-k}\Omega^k S^n\{[2]\}$, which projects to a class by the same
name in $H_*\Omega^k S^n$. There is a unique non-zero class
$x_{n-k-1}$ in $H_{n-k-1}\Omega^k S^n\{[2]\}$ which is in the image
of the natural map
$$H_*\Omega^{k+1}S^n \to\ H_*\Omega^k S^n\{[2]\}$$ such that
$Sq^1_*(x_{n-k}) = x_{n-k-1}$. The Nishida relations are given by
$$Sq^t_*Q_r(x) = \Sigma_{0 \leq i}(t-2i,r+q-2t+2i)Q_{r-t+2i}Sq^i_*(x)$$ for any class $x$ of degree
$q$ with binomial coefficients given by $$(a,b) =(a+b)!/a! \cdot
b!$$ for $a,b \geq 0$. The Steenrod operations in $H_*\Omega^k
S^n\{[2]\}$ then follow by specialization. Examples of Steenrod
operations acting on the classes
\begin{enumerate}
  \item $Q_i(x_{n-k})$ for $0 \leq i \leq k-1$ and
  \item $Q_j(x_{n-k-1})$ for $0 \leq j \leq k$
\end{enumerate} will be considered next. Specialize to the case $H_*\Omega^{2^n} S^{q2^{n+2} + 1}\{[2]\}$,
$q \geq 1$. Thus there are classes $$v = x_{q2^{n+2}-2^n+1}$$ and
$$u = x_{q2^{n+2}-2^n}$$ with $Sq^1_*x_{q2^{n+2}-2^n+1}  = Sq^1_*(v)
= u = x_{q2^{n+2}-2^n}$ as given by the above remarks ( for which
reference to degrees is deliberately omitted in the cases of $u$,
and $v$ ). Notice that

\begin{enumerate}
        \item $Sq^{2^{n}}_* Q_{2^n-1}(v)
        = \Sigma_{0 \leq i}(2^{n}-2i, 2^n -1 + |v| -2^{n+1}+2i)Q_{-1+2i}
        Sq^i_*v$,
        \item $Sq^{2^{n}}_* Q_{2^n-1}(v)
        = (2^{n}-2, 2^n - 1 + |v| -2^{n+1}+2)Q_{1}Sq^1_*v$ and
        \item $Sq^{2^{n}}_* Q_{2^n-1}(v)
        = (2^{n}-2,q2^{n+2}-2^{n+1}+2)Q_{1}u$.
\end{enumerate}

By \cite{CCPS}, the module of primitive elements in $H_{*}\Omega^k
S^n\{[2]\}$ is spanned by $Q_{i_1}Q_{i_2} \cdots Q_{i_k}(v)$ and
$Q_{j_1}Q_{j_2} \cdots Q_{i_m}(u)$ for $1 \leq i_1 \leq i_2 \leq
\cdots \leq i_k \leq n-1$ and $1 \leq j_1 \leq j_2 \leq \cdots \leq
j_m \leq n$ as well as their $2^r$-th. The degree of $Q_i(x)$,
$|Q_i(x)|$, is given by $i + 2|x|$. If $2 < 2k < n-2$, then the next
result follows by induction and a standard degree count \cite{ CCPS,
1286}

\begin{lem}\label{lemma:primitives}
Assume that $$2 < 2k < n-2.$$ There are unique non-trivial primitive
elements $$v \epsilon H_{n-k}\Omega^k S^n\{[2]\}$$ and $$u \epsilon
H_{n-k-1}\Omega^k S^n\{[2]\}.$$ A basis for the module of primitives
$PH_*\Omega^k S^n\{[2]\}$ in degrees less than $3(n-k-1)$ is
$$\{u, v, Q_{i}(u), Q_{j}(v)| 0 \leq i \leq n-k, 0 \leq j
\leq n-k-1\}.$$ Furthermore, the element $Q_{k-1}(v)$ is the unique
non-trivial primitive in $H_{2n-k-1}\Omega^k S^n\{[2]\}$. Thus there
is exactly one non-trivial primitive element in $H_{q2^{n+3}-2^n
-1}\Omega^{2^n} S^{q2^{n+2} + 1}\{[2]\}$ given by $Q_{2^n-1}(v)$.
\end{lem}

The action of certain Steenrod operations is given next.
\begin{lem}\label{lemma:squares on fibre of degree 2}
If $n>1$ and $q\geq 1$, then $(2^n-2, q2^{n+2}-2^{n+1}+2)=0$ modulo
$2$. Thus in $H_*\Omega^{2^n} S^{q2^{n+2} + 1}\{[2]\}$, there are
unique non-trivial primitive elements $v = x_{q2^{n+2}-2^n+1}$ and
$u = x_{q2^{n+2}-2^n}$ with $Sq^1_*v = u$ and $Sq^{2^{n}}_*
Q_{2^n-1}(v) = 0$ for the unique non-zero primitive element
$Q_{2^n-1}(v)$ in degree $q2^{n+3}-2^n-1$.
\end{lem}

\begin{proof}
To prove Lemma \ref{lemma:squares on fibre of degree 2}, recall the
well-known method for evaluating binomial coefficients via $p$-adic
expansions \cite{Steenrod.Epstein}: Let $p$ be a prime. Assume that
$a$ and $b$ are strictly positive integers for which there are
choices of $p$-adic expansions $a=\sum_{i=0}^m a_i p^i$ and
$b=\sum_{i=0}^m b_i p^i$, $0\leq a_i, b_i < p$. Then $$( {}^a_b ) =
\prod_{i=0}^m ( {}^{a_i}_{b_i} )\pmod p. $$

Thus, consider the mod $2$ reduction of the binomial coefficient
given by $$(2^n-2, q2^{n+2}-2^{n+1}+2)= \binom \alpha \beta$$ for
$\alpha =  q2^{n+2}-2^n$ and $\beta = 2^n-2$. Notice that
$$q2^{n+2}-2^n = (q-1)2^{n+2}+ 2^{n+2}-2^{n} = (q-1)2^{n+2}+
2^{n+1}+2^{n}$$ where $q-1$ is a non-negative integer. In this case,
the $2$-adic expansion for $q2^{n+2}-2^n =\sum_{i=0}^m a_i 2^i$ has
$a_i=0$  for $i\leq n-1$ and $a_n=a_{n+1}=1$.  The $2$-adic
expansion for $2^n-2 = \sum_{i=0}^m b_i 2^i$ has $b_i=1$ for $1\leq
i \leq n-1$ and all other $b_i$ are 0. Hence, it follows that $1
\leq i \leq n-1$
$$\left({}^{a_i}_{b_i}\right) = \left({}^0_1\right)\equiv 0 \pmod
2,$$ and thus $$( {}^{2^{n+2}q-2^n}_{2^n-2} ) = \prod_{i=0}^m (
{}^{a_i}_{b_i} )=0\pmod 2.$$

\end{proof}

Features of the homology of $\Omega^k S^n\{\Psi\}$ were worked out
in \cite{CCPS} for $ 1 < k < n-3 $ using the stabilization map
$E:S^n \to\ QS^n$ to obtain a map $$\gamma:\Omega^k S^n\{\Psi\} \to\
(\Omega^k QS^n)\{\Psi\}.$$ Computations with the Steenrod operations
will be given using $\gamma$. The following properties are satisfied
in these cases by \cite{CCPS}.

\begin{enumerate}
  \item The mod-$2$ homology of $\Omega^k S^n\{\Psi\}$
is isomorphic to $$H_*\Omega^k S^n \otimes H_*\Omega^{k+1}S^n$$ as a
Hopf algebra. The natural map $$H_*\Omega^k S^n\{\Psi\} \to
H_*\Omega^k S^n$$ induced by a $k-1$-fold loop map and is an
epimorphism of Hopf algebras over the mod-$2$ Steenrod algebra. Thus
there are unique (non-trivial) primitive elements $x_{j}$ in
$H_{j}\Omega^k S^n\{\Psi\}$ for $j$ equal to either $n-k$, or $n-k
-1$.

It is convenient for the computations below to abbreviate the names
of two elements in the following way: $u= x_{n-k-1}$ and $v=
x_{n-k}$.

\item The mod-$2$ homology of $(\Omega^k QS^n)\{\Psi\}$
is isomorphic to $$H_*\Omega^k QS^n \otimes H_*\Omega^{k+1}QS^n$$ as
a Hopf algebra. The natural map $$H_* (\Omega^k QS^n)\{\Psi\} \to
H_* \Omega^k QS^n$$ induced by an infinite loop map is an
epimorphism of Hopf algebras over the mod-$2$ Steenrod algebra. Thus
there are unique (non-trivial) primitive elements $y_{j}$  in
$H_j((\Omega^k QS^n)\{\Psi\})$ for $j$ equal to either $n-k$, or
$n-k -1$.

It is again convenient for the computations below to abbreviate the
names of two elements in the following way: $w= y_{n-k-1}$ and $z=
y_{n-k}$.

\item There exist primitive elements in $H_{i+
2(n-k)}\Omega^k S^n\{\Psi\}$ for $0 \leq i \leq k-1$ denoted $\bar
Q_i(x_{n-k})$ which project to elements $Q_{i}(x_{n-k})$ in
$H_*(\Omega^k  S^n)$.

The elements $\bar Q_i(x_{n-k})$, for $0 \leq i \leq k-2$ are given
by the Araki-Kudo-Dyer-Lashof operations $Q_i$ on $x_{n-k}$.
However, the element $\bar Q_{k-1}(x_{n-k})$ is not given by an
operation. The symbol $\bar Q_{k-1}(-)$ is a formal bookkeeping
device; this symbol does not mean that it is given by an operation.

\item By \cite{CCPS}, the formula  $$\gamma_*(\bar
Q_{k-1}(x_{n-k})) = Q_{k-1}(y_{n-k}) +Q_{k+1}(y_{n-k+1})$$ is
satisfied in case $n$ is not equal to $3$ modulo $4$.
\end{enumerate}

Using this information, the Steenrod operations on $\bar
Q_{k-1}(x_{n-k})$ in $H_*\Omega^k S^n\{\Psi\}$ will be given from
the Nishida relations by using the map $\gamma$ in the special cases
$\Omega^{2^n} S^{q2^{n+2} + 1}\{\Psi\}$ with $n \geq 1$. Thus
$q2^{n+2} + 1$ is not $3$ modulo $4$, and \cite{CCPS} applies. A
direct count of degrees analogous to that in Lemma
\ref{lemma:primitives} gives uniqueness of certain primitives.

\begin{lem}\label{lemma:primitives.two}
Assume that $$2 < 2k < n-2.$$ A basis for the module of primitives
$PH_*\Omega^k S^n\{\Psi\}$ in degrees less than $3(n-k-1)$ is
$$\{x_{n-k-1}, x_{n-k}, Q_{i}(x_{n-k-1}), \bar Q_{j}(x_{n-k})| 0 \leq i \leq n-k, 0 \leq j
\leq n-k-1\}.$$ Furthermore, the element $\bar Q_{k-1}(x_{n-k})$ is
the unique non-trivial primitive of degree $2n-k-1$.  Thus there is
exactly one non-trivial primitive element in $$H_{q2^{n+3}-2^n
+1}\Omega^{2^n} S^{q2^{n+2} + 1}\{\Psi\}$$ given by $\bar
Q_{2^n-1}(x_{n-k})$.
\end{lem}

\begin{lem}\label{lemma:squares on fibres}
If $n>1$ and $q\geq 1$, then $(2^n, q2^{n+2}-2^{n+1}+1) =1 $ modulo
$2$. Hence, the unique non-zero primitive elements $u =
x_{q2^{n+2}-2^n}$, and $v = x_{q2^{n+2}-2^n+1}$ in $H_*\Omega^{2^n}
S^{q2^{n+2} + 1}\{\Psi\}$ satisfy
\begin{enumerate}
  \item $Sq^1_*v  = u$,
  \item $Sq^{2^{n}}_* \bar Q_{2^n-1}(v) \neq 0$ and
  \item $\bar Q_{2^n-1}(v)$ is the unique non-zero primitive element
  in $$H_{q2^{n+3}-2^n-1}(\Omega^{2^n} S^{q2^{n+2} + 1}\{\Psi\}).$$
\end{enumerate}
\end{lem}

\begin{proof}[Proof of \ref{lemma:squares on fibres}]

To prove Lemma \ref{lemma:squares on fibres}, first consider the
binomial coefficient $$(2^n, q2^{n+2}-2^{n+1}+1)= \left(
{}^{q2^{n+2}-2^{n+1}+2^n+1}_{2^n} \right).$$ Assuming $q\geq 1$ the
2-adic expansion for $q2^{n+2}-2^{n+1}+2^n+1= (q-1)2^{n+2} + 2^{n+1}
+ 2^n +1 = \sum_{i=0}^m a_i 2^i$ has $a_{n+1}=a_n=a_0=1$ with $a_i =
0$ for $1 \leq i \leq n-1$.  The 2-adic expansion for
$2^n=\sum_{i=0}^m b_i 2^i$ has $b_n=1$ with all other $b_i = 0$.
Thus if $ i \neq n$, $\left( {}^{a_i}_{b_i} \right)=
\left({}^{a_i}_0\right)=1$ with $\left( {}^{a_n}_{b_n}
\right)=\left({}^1_1\right)$.  Thus $\left(
{}^{2^{n+2}q-2^{n+1}+2^n+1}_{2^n} \right) = \prod_{i=0}^m (
{}^{a_i}_{b_i} ) = 1 \pmod 2$ and the formula for binomial
coefficients follows.

Recall the abbreviation of the names of classes as above with $u=
x_{n-k-1}$ in $H_{n-k-1}\Omega^k S^n\{\Psi\}$ and $v= x_{n-k}$ in
$H_{n-k}\Omega^k S^n\{\Psi\}$. A second abbreviation is given by $w
= y_{n-k-1}$ in $H_{n-k-1}((\Omega^k QS^n)\{\Psi\})$, and $z =
y_{n-k}$ in $H_{n-k}((\Omega^k QS^n)\{\Psi\})$.

Next consider the element $$\gamma_*(\bar Q_{2^n-1}(v))=
Q_{2^n-1}(z)+ Q_{2^n+1}(w).$$ The next properties follow at once.
\begin{enumerate}
        \item $Sq^{2^{n}}_* \gamma_*(\bar Q_{2^n-1}(v)) =
        Sq^{2^{n}}_*Q_{2^n-1}(z)
        +Sq^{2^{n}}_*Q_{2^n+1}(w)$ for $n-k = q2^{n+2}-2^n$.
        \item  $Sq^{2^{n}}_*(Q_{2^n-1}(z))= 0$ by Lemma
        \ref{lemma:squares on fibre of degree 2}.
    \item $Sq^{2^{n}}_* Q_{2^n+1}(w) = (2^{n},2^n+1 +
    q2^{n+2}-2^n - 2^{n+1})Q_{2^n+1- 2^n}(w)$.
    \item Since the binomial coefficient $(2^n, q2^{n+2}-2^{n+1}+1)$
is $1$ modulo two, $$Sq^{2^{n}}_* Q_{2^n+1}(w) = Q_{1}(w)$$ and so
$Sq^{2^{n}}_* Q_{2^n+1}(w) \neq 0$.

\item The element $Sq^{2^{n}}_* \gamma_*(\bar Q_{2^n-1}(v))$ is
non-zero.

\item By Lemma \ref{lemma:primitives.two}, the element $\bar Q_{2^n-1}(v)$ is the
unique non-trivial primitive in $H_{q2^{n+3}-2^n-1}(\Omega^{2^n}
S^{q2^{n+2} + 1}\{\Psi\})$.
\end{enumerate}
 The lemma follows.
\end{proof}

The proof of Theorem \ref{thm:theorem two} is given next.
\begin{proof}[Proof of \ref{thm:theorem two}]

Lemmas \ref{lemma:squares on fibre of degree 2}, and
\ref{lemma:squares on fibres} immediately imply Theorem
\ref{thm:theorem two} that $\Omega^{2^n} S^{2^{n+2}q + 1}\{2\}$, and
$\Omega^{2^n} S^{2^{n+2}q + 1}\{[2]\}$ are not homotopy equivalent
for $n > 1$ and $q \geq 1$  as these spaces have different actions
of the Steenrod algebra as follows.

Assume that $n>1$ and $q\geq 1$. By \ref{lemma:squares on fibre of
degree 2}, there is an unique non-zero primitive element in
$H_t\Omega^{2^n} S^{q2^{n+2} + 1}\{[2]\}$ for $t = q2^{n+3}-2^n-1$
given by $Q_{2^n-1}(v)$. Furthermore, this element satisfies
$Sq^{2^{n}}_* Q_{2^n-1}(v) = 0$.

Again assume that $n>1$ and $q\geq 1$. By \ref{lemma:squares on
fibres} there is an unique non-zero primitive element in
$H_t\Omega^{2^n} S^{q2^{n+2} + 1}\{\Psi\}$ for $t = q2^{n+3}-2^n-1$
given by $\bar Q_{2^n-1}(z)$. Furthermore, this element satisfies
$Sq^{2^{n}}_* \bar Q_{2^n-1}(z) \neq 0$.

The unique non-zero primitive elements in degree $q2^{n+3}-2^n-1$
support different actions of $Sq^{2^{n}}_*$. Hence the mod-$2$
cohomology of the spaces $\Omega^{2^n} S^{q2^{n+2} + 1}\{2\}$, and
$\Omega^{2^n} S^{q2^{n+2} + 1}\{\Psi\}$ differ. Theorem
\ref{thm:theorem two} follows

\end{proof}

\section{On the Proof of Proposition \ref{prop:prop four}}

A proof of the main part of Proposition \ref{prop:prop four}
\cite{1286} is given below for convenience of the reader.
Proposition \ref{prop:prop four} is a restated version Proposition
$11.3$ of \cite{1286} in which there is a misprint where
$\Omega^q(\phi)$ should be $\Omega^{q-1}(\phi)$ ( as stated in
section $1$ here ).
\begin{thm}\label{thm:null-homotopies}
Assume that the composite
\[
\begin{CD}
\Omega^kS^{2n+1}@>{\Omega^{k-1}h_2}>>
\Omega^kS^{4n+1}@>{\Omega^k(w_{2n+1})}
>>\Omega^kS^{2n+1}
\end{CD}
\] is null-homotopic.

Then the composite
\[
\begin{CD}
\Sigma^{2n+1-k}(\mathbb R\mathbb P^{2n}_{2n-k+1})@>{collapse}>>
S^{4n+1} @>{w_{2n+1}}>> S^{2n+1}
\end{CD}
\] is null-homotopic.
\end{thm}

{\bf Remark: \label{remark:null-homotopies}} The results in this
article do not rule out the possibility that the converse of Theorem
\ref{thm:null-homotopies} may be satisfied.

\begin{proof}[Proof of Theorem \ref{thm:null-homotopies}]
Let $F_2 = F_2(\Omega^i S^{2n+1})$ denote the second filtration of
the May-Milgram construction for $\Omega^k S^{2n+1}$ with $I:F_2
\to\ \Omega^k S^{2n+1}$ giving the natural inclusion. One fact is
that there is a cofibration sequence
\[
\begin{CD}
S^{2n+1-k} @>{inclusion}>>  F_2 @>{collapse}>> \Sigma^{2n+1-k}
(\mathbb R\mathbb P^{2n}_{2n-k+1})@>{\delta}>> S^{2n+2-k}@>{}>>
\cdots
\end{CD}
\] with the property that the $k$-fold suspension of $\delta$,
$\Sigma^k(\delta)$, is null-homotopic.

Next, consider the commutative diagram

\[
\begin{CD}
S^{2n+1-k} @>{inclusion}>> F_2 \\
 @V{identity}VV         @VV{I}V     \\
 S^{2n+1-k}       @>{E^k}>> \Omega^kS^{2n+1} \\
 @V{identity}VV         @VV{\Omega^{k-1}h_2}V     \\
 S^{2n+1-k}       @>{}>> \Omega^kS^{4n+1}. \\
\end{CD}
\] Notice that if $n \geq 1$, then any map
$$S^{2n+1-k} \to\  \Omega^kS^{4n+1}$$ is null and so there
is a homotopy commutative diagram

\[
\begin{CD}
 F_2 @>{collapse}>> \Sigma^{2n+1-k} (\mathbb R\mathbb P^{2n}_{2n-k+1}) \\
 @V{\Omega^{k-1}h_2\circ I}VV         @VV{\Theta}V     \\
\Omega^kS^{4n+1}     @>{identity}>> \Omega^kS^{4n+1} \\
\end{CD}
\]for some map $\Theta$ by the standard properties of the
cofibration sequence above for $F_2$. Observe that there is a
homotopy commutative diagram

\[
\begin{CD}
\Sigma^{2n+1-k} (\mathbb R\mathbb P^{2n}_{2n-k+1}) @>{1}>>
\Sigma^{2n+1-k} (\mathbb R\mathbb P^{2n}_{2n-k+1}) \\
 @V{K}VV         @VV{\Theta}V     \\
S^{4n+1-k}     @>{E^k}>> \Omega^kS^{4n+1} \\
\end{CD}
\] as $\Omega^kS^{4n+1}$ is $4n-k$-connected.

Next assume that $$\Omega^k(w_{2n+1})\circ \Omega^{k-1}h_2$$ is
null. This assumption gives that $\Omega^k(w_{2n+1})\circ
\Omega^{k-1}h_2 \circ I$ is also null-homotopic. Thus there is yet
another homotopy commutative diagram
\[
\begin{CD}
 F_2 @>{collapse}>> \Sigma^{2n+1-k} (\mathbb R\mathbb P^{2n}_{2n-k+1}) \\
 @V{\Omega^{k-1}h_2\circ I}VV         @VV{\Theta}V     \\
\Omega^kS^{4n+1}     @>{identity}>> \Omega^kS^{4n+1} \\
@V{\Omega^k(w_{2n+1})}VV  @VV{\Omega^k(w_{2n+1})}V     \\
\Omega^kS^{2n+1}     @>{identity}>> \Omega^kS^{2n+1} \\
\end{CD}
\] where the vertical left-hand composite is null by assumption.
Thus the right-hand vertical composite map in this diagram

\[
\begin{CD}
\Sigma^{2n+1-k}(\mathbb R\mathbb P^{2n}_{2n-k+1}) @>{\Theta}>>
\Omega^kS^{4n+1} @>{\Omega^k(w_{2n+1})}>> \Omega^kS^{2n+1}
\end{CD}
\]factors, up to homotopy, through the cofibre of
the natural map
$$F_2 \to\ \Sigma^{2n+1-k} (\mathbb R\mathbb P^{2n}_{2n-k+1}).$$
Hence there is a homotopy commutative diagram

\[
\begin{CD}
\Sigma^{2n+1-k} (\mathbb R\mathbb P^{2n}_{2n-k+1}) @>{\delta}>> S^{2n+2-k} \\
@VV{\Omega^k(w_{2n+1}) \circ \Theta}V  @VV{\alpha}V     \\
\Omega^kS^{2n+1}     @>{identity}>> \Omega^kS^{2n+1} \\
\end{CD}
\] for some choice of map $\alpha$.

Since the cofibration sequence
\[
\begin{CD}
S^{2n+1-k} @>{inclusion}>> F_2 @>{collapse}>> \Sigma^{2n+1-k}
(\mathbb R\mathbb P^{2n}_{2n-k+1})
 @>{\delta}>> S^{2n+2-k}
\end{CD}
\] satisfies the property that $\Sigma^k(\delta)$ is null, passage to adjoints
gives the following homotopy commutative diagram.

\[
\begin{CD}
\Sigma^k(\Sigma^{2n+1-k} \mathbb R\mathbb P^{2n}_{2n-k+1})
@>{\Sigma^k(\delta)}>>
\Sigma^k(S^{2n+2-k}) \\
@V{\Sigma^k(\Omega^k(w_{2n+1}) \circ \Theta)}VV
@VV{\Sigma^k(\alpha)}V     \\
\Sigma^k(\Omega^kS^{2n+1})     @>{1}>> \Sigma^k(\Omega^kS^{2n+1})\\
@V{evaluation}VV        @VV{evaluation}V     \\
S^{2n+1}   @>{1}>> S^{2n+1}.  \\
\end{CD}
\] Hence the vertical left-hand composite
factors through the null-homotopic map $\Sigma^k(\delta)$, and the
theorem follows.

\end{proof}

\bibliographystyle{amsalpha}

\begin{thebibliography}{99}

\bibitem{BP} E.~H.~Brown, and F.~P.~Peterson, {\em Whitehead products, and cohomology},
Quart. J. Math. {\bf15 }(1964), 116-120.

\bibitem{CCPS} H.~E.~A.~Campbell, F.~R.~Cohen, F.~P.~Peterson,
and P.~S.~Selick, {\em The space of maps of Moore spaces to
spheres}, Annals of Math. Studies, \textbf{13}(1987), 72--100.

\bibitem{1286}F.~R.~Cohen, {\em A course in some aspects of
classical homotopy theory}, S.L.M., \textbf{1286}(1987), 1--92.

\bibitem{Adem proceedings}F.~R.~Cohen, {\em On the Whitehead square, Cayley-Dickson algebras, and rational functions},
Boletin de la Sociedad Mexicana, \textbf{37}(1992), 55-62.

\bibitem{Harper} J.~Harper, {\em Secondary cohomology operations },  Graduate Studies in Mathematics, Amer. Math.
Soc., {\bf 49}(2002).

\bibitem{Hunter} T.~Hunter, {\em On $H\sb *(\Omega^{n+2}S^{n+1}; F\sb 2)$},
Trans. Amer. Math. Soc., {\bf 314}(1989), no. 1, 405--420.

\bibitem{JM} I.~Johnson, J.~Merzel, {\em A Class of Left Ideals of the Steenrod Algebra}, preprint.

\bibitem{Kochman} S.~O.~Kochman, {\em Stable homotopy groups of spheres}, S.L.M., {\bf 1423}(1990).

\bibitem{Mahowald} M.~Mahowald, {\em Some remarks on the Kervaire invariant form the homotopy point of view},
Proc. of Symposia in Pure Math. {\bf 22}(1971), Providence R. I.,
165-169.

\bibitem{May} J.~P.~May, {\em The Geometry of Iterated Loop Spaces}, S.L.M., \textbf{268}(1972).

\bibitem{Moreno} G.~Moreno, {\em The zero divisors of the Cayley-Dickson algebras over the real numbers}, preprint.

\bibitem{Nishida} G.~Nishida, {\em Cohomology operations in iterated loop spaces }, Proc. Japan Acad., \textbf{44}(1968),
104-109.

\bibitem{Selick} P.~S.~Selick, {\em A reformulation of the Arf invariant one mod p problem and applications to atomic spaces},
Pac. J. Math., \textbf{108}(1983), 431-450.

\bibitem{Steenrod.Epstein} N.~E.~Steenrod, and D.B.~Epstein, {\em Cohomology Operations}, Annals of Math. Studies,
 {\bf 50}(1968).

 \bibitem{Toda} H.~Toda, {\em Composition Methods in the Homotopy Groups of Spheres}, Annals of Math. Studies,
 \textbf{49}(1968).

\end{thebibliography}

\end{document}